\documentclass[11pt,reqno]{amsart}
\usepackage[T1]{fontenc}
\usepackage{lmodern,microtype}
\usepackage{enumitem}
\usepackage{amsmath,amssymb,amsfonts,amsthm,mathtools,mathrsfs}
\usepackage[margin=1.08in]{geometry}
\usepackage{xcolor}
\usepackage[colorlinks=true,linkcolor=blue!45!black,citecolor=blue!45!black,urlcolor=blue!45!black,pdftitle={Two Routes from Total Curvature to Topology of Complete Euclidean Immersions}]{hyperref}
\usepackage[nameinlink,capitalize,noabbrev]{cleveref}
\newtheorem{theorem}{Theorem}[section]
\newtheorem{proposition}[theorem]{Proposition}
\newtheorem{lemma}[theorem]{Lemma}
\newtheorem{corollary}[theorem]{Corollary}
\theoremstyle{definition}
\newtheorem{definition}[theorem]{Definition}
\theoremstyle{remark}
\newtheorem{remark}[theorem]{Remark}

\newcommand{\R}{\mathbb R}
\newcommand{\C}{\mathbb C}
\newcommand{\F}{\mathbb F}
\newcommand{\cH}{\mathcal H}
\newcommand{\TT}{\mathcal T}
\newcommand{\SSS}{\mathbb S}
\newcommand{\TA}{\operatorname{TA}}

\newcommand{\rank}{\operatorname{rank}}

\newcommand{\bb}{\mathbf b}
\newcommand{\np}{\mathrm{np}}
\DeclareMathOperator{\supp}{supp}
\title[Heat kernels and Morse theory for Euclidean immersions]{Total Curvature and Topology of Complete Euclidean Immersions: Heat Kernel and Boundary Morse Approaches}
\author{Lei Zhang}
\address{Shanghai Center for Mathematical Sciences, Fudan University, Shanghai 200438, China}

\email{22110840013@m.fudan.edu.cn}

\begin{document}
\begin{abstract}
Two mechanisms relate the scale-invariant integral of the second fundamental form to topology at infinity of a complete Euclidean immersion. The first retains the mean-curvature potential in the Michael--Simon Sobolev inequality, establishes an ultracontractive estimate for the semigroup generated by $\Delta+|H|^2$, and uses a vector-bundle Cwikel--Lieb--Rozenblum bound to estimate spaces of $L^2$-harmonic forms. Together with the theory of nonparabolic ends, this yields $e(M)\le 1+C_n\int_M|A|^n$ for arbitrary complete immersed $n$-manifolds, $n\ge3$. The second method, under the additional assumption of finitely many single-sheeted $C^2$-asymptotically planar graphical ends, uses a boundary Morse complex and the normal Gauss-map area formula to prove $\sum_j b_j(M;\F)+e(M)-2\le\TA(F)$. It supplies an explicit end coefficient and bounds ordinary Betti numbers. We provide proofs of both mechanisms, identify the exact role of their respective hypotheses, and discuss minimal and Hamiltonian stationary Lagrangian specializations.
\end{abstract}
\maketitle

\section{Introduction and main statements}\label{sec:intro}
Let $F:M^n\to\R^N$, $N>n\ge3$, be a smooth complete connected isometric immersion without boundary. Unless otherwise stated, the immersion need not be minimal, proper, embedded, or orientable. Let $A$ be the second fundamental form and adopt the trace convention
\[
 H=\sum_{i=1}^n A(e_i,e_i),\qquad
 \TT(F):=\int_M|A|^n\,d\mu.
\]
The latter quantity is invariant under Euclidean dilations. We write $e(M)$ for the number of Freudenthal ends. For $0\le p\le n$ let
\[
 \cH^p_{(2)}(M)=\{\omega\in L^2(\Lambda^pT^*M):d\omega=d^*\omega=0\},\qquad
 b^p_{(2)}(M)=\dim\cH^p_{(2)}(M).
\]
Completeness identifies these spaces with reduced $L^2$-cohomology. They should not be confused with the ordinary Betti numbers
\[
 \bb(M;\F)=\sum_{p=0}^n\dim_{\F}H_p(M;\F)
\]
over a coefficient field $\F$.

The central distinction is between an analytic conclusion available under finite $\TT(F)$ alone and an ordinary-topology conclusion which requires geometric information about the ends.

\begin{theorem}[Heat-kernel method: general immersion]\label{thm:general}
Let $F:M^n\to\R^N$ be as above and assume $\TT(F)<\infty$. Every end of $M$ is nonparabolic. There exist finite constants $C_n$ and $C_{n,p}$, depending only on the indicated intrinsic dimensions, such that
\begin{align}
 b^p_{(2)}(M)&\le C_{n,p}\TT(F),\qquad 1\le p\le n-1,\label{eq:lpbound}\\
 e(M)&\le1+C_n\TT(F).\label{eq:generalend}
\end{align}
In particular $M$ has finitely many ends. If $\TT(F)$ is sufficiently small, $M$ has exactly one end and all intermediate reduced $L^2$-Betti numbers vanish.
\end{theorem}

The first method exploits the mean-curvature term in the Michael--Simon inequality rather than trying to absorb it. The resulting Schr\"odinger heat kernel admits a dimensional $t^{-n/2}$ upper bound without any smallness assumption on $\|H\|_{L^n}$. A vector-bundle version of the abstract CLR estimate then controls the Hodge kernels; the end conclusion uses the Fu--Li nonparabolicity criterion and the Li--Tam construction. The shifted spectral estimate is formulated explicitly below, since it is not literally the ordinary-Laplacian formulation of B\'erard--Besson's Theorem 37 \cite{BerardBesson}.

Put $k=N-n$ and define
\begin{equation}\label{eq:cconst}
 c_{n,N}:=\frac{|\SSS^{k-1}|}{n^{n/2}|\SSS^{N-1}|}.
\end{equation}
For the second method, we use the graphical-end hypothesis defined precisely in \cref{def:ends}.

\begin{theorem}[Morse-complex method: planar graphical ends]\label{thm:ordinary}
Let $F:M^n\to\R^N$ be a complete connected smooth immersion with finitely many single-sheeted $C^2$-asymptotically planar graphical ends in the sense of \cref{def:ends}. Then for every field $\F$,
\begin{equation}\label{eq:morsecl}
 \bb(M;\F)+e(M)-2\le\TA(F),
\end{equation}
where $\TA(F)$ denotes the normalized total absolute curvature. In particular,
\[
 2(e(M)-1)\le\TA(F).
\]
If $\TT(F)<\infty$, then
\begin{align}
 \bb(M;\F)&\le1+c_{n,N}\TT(F),\label{eq:ordinarybound}\\
 e(M)&\le1+\tfrac12c_{n,N}\TT(F).\label{eq:planarend}
\end{align}
\end{theorem}

\begin{corollary}[Minimal immersions]\label{cor:minimal}
For a complete connected minimal immersion $F:M^n\to\R^N$ with $N>n\ge3$ and $\TT(F)<\infty$, the conclusions of both theorems hold. The immersion is proper, the manifold has finite topological type, and
\[
 e(M)\le1+\tfrac12c_{n,N}\TT(F),\qquad
 \bb(M;\F)\le1+c_{n,N}\TT(F).
\]
\end{corollary}

Here the geometric input is the  regular-at-infinity theorem for minimal immersions in arbitrary codimension \cite{Anderson,DingZhang}, not merely the formal identity $H=0$.

\begin{corollary}[Hamiltonian stationary Lagrangians]\label{cor:hsl}
Let $L^n\to\C^n\cong\R^{2n}$ be a complete connected Hamiltonian stationary Lagrangian immersion, $n\ge3$, satisfying $\int_L|A|^n<\infty$. The general end and intermediate reduced $L^2$-cohomology bounds in \cref{thm:general} hold. If, additionally, $L$ has the graphical ends of \cref{def:ends}, then the ordinary Betti-number and explicit end estimates of \cref{thm:ordinary} hold with $N=2n$.
\end{corollary}

The identities $H=J\nabla\beta$ and $\Delta_L\beta=0$ hold locally for a Hamiltonian stationary Lagrangian, but are not needed by the general spectral argument. They do not on their own furnish the required $C^2$ graphical-end structure.

\subsection*{Context and prior work}
The potential-theoretic and spectral ingredients originate in Michael--Simon \cite{MichaelSimon}, B\'erard--Besson \cite{BerardBesson}, Li--Tam \cite{LiTam}, and the noncompact $L^2$-cohomology literature \cite{Carron,FuLi,FuXu,CMV}. Quantitative harmonic-index bounds for minimal hypersurfaces were obtained by Mei--Xu \cite{MeiXu}. The boundary Morse construction is that of Laudenbach \cite{Laudenbach}; the integral-geometric framework traces to Chern--Lashof \cite{ChernLashof}. The minimal-end input is available in Anderson's work and the recent Ding--Zhang formulation \cite{Anderson,DingZhang}. We compare the mechanisms and their limitations in \cref{sec:comparison}, without assuming any general equivalence between ordinary and reduced $L^2$-cohomology.
\subsection{{Statement on AI usage}}  The author used  ChatGPT as an assistive tool in the development of the paper. The author takes full responsibility for all mathematical arguments.
\section{The mean-curvature Schr\"odinger heat kernel}
The starting point is the Euclidean Michael--Simon inequality \cite{MichaelSimon}, in which no minimality hypothesis is needed.

\begin{lemma}[Sobolev inequality with a mean-curvature term]\label{lem:MS}
There is $S_n<\infty$ such that for all $u\in C_c^1(M)$,
\begin{equation}\label{eq:MS}
 \|u\|_{L^{2^*}(M)}^2\le S_n\int_M\left(|\nabla u|^2+|H|^2u^2\right)\,d\mu,
 \qquad 2^*=\frac{2n}{n-2}.
\end{equation}
\end{lemma}
\begin{proof}
The Michael--Simon inequality reads $\|v\|_{n/(n-1)}\le C_n\int(|\nabla v|+|H|\,|v|)$ for nonnegative compactly supported $v$. Apply it to $v=|u|^\gamma$, $\gamma=2(n-1)/(n-2)$, first with a harmless smooth approximation at zeros. H\"older's inequality yields
\[
 \|u\|_{2^*}^{\gamma}\le C_n\|u\|_{2^*}^{\gamma-1}\bigl(\|\nabla u\|_2+\||H|u\|_2\bigr).
\]
After division and squaring, we obtain \eqref{eq:MS}.
\end{proof}

Set $q=|H|^2$ and let $P_q=\Delta+q$ be the nonnegative Friedrichs operator associated with
\[
 Q_q[u]=\int_M\bigl(|\nabla u|^2+qu^2\bigr)\,d\mu.
\]
On a complete manifold, the form closure of $C_c^\infty(M)$ is the natural realization considered here.

\begin{proposition}[Dimensionally uniform heat-kernel estimate]\label{prop:heat}
The minimal scalar heat kernel $p_q(t,x,y)$ of $e^{-tP_q}$ satisfies
\begin{equation}\label{eq:heat}
 0\le p_q(t,x,x)\le a_n t^{-n/2},\qquad t>0,\ x\in M,
\end{equation}
with $a_n$ independent of the immersion, $N$, and $\int_M|H|^n$.
\end{proposition}
\begin{proof}
Interpolation followed by \eqref{eq:MS} gives the Nash inequality
\begin{equation}\label{eq:Nash}
 \|u\|_2^{2+4/n}\le \|u\|_1^{4/n}\|u\|_{2^*}^2
 \le S_nQ_q[u] \|u\|_1^{4/n}.
\end{equation}
As $q\ge0$, the Schr\"odinger semigroup is sub-Markovian: $\|e^{-tP_q}f\|_1\le\|f\|_1$. If $Y(t)=\|e^{-tP_q}f\|_2^2$, then
\[
 Y'(t)=-2Q_q[e^{-tP_q}f]\le-2S_n^{-1}\|f\|_1^{-4/n}Y(t)^{1+2/n}.
\]
Integration implies $\|e^{-tP_q}\|_{1\to2}\le C_nt^{-n/4}$. Self-adjointness and the semigroup identity give $\|e^{-tP_q}\|_{1\to\infty}\le a_nt^{-n/2}$. This implies \eqref{eq:heat}, using a continuous representative of the heat kernel. The argument is the usual Nash--ultracontractivity implication \cite{Davies}.
\end{proof}

\section{A bound-state estimate with a reference potential}
For clarity we spell out the extension of the noncompact B\'erard--Besson inequality that is used here. This is an abstract heat-kernel bound-state argument, not an application of their theorem to the ordinary heat kernel of $M$.

\begin{proposition}[Shifted CLR estimate for bundles]\label{prop:CLR}
Suppose $n\ge3$, $q\ge0$, and $p_q(t,x,x)\le at^{-n/2}$ on a complete $n$-manifold. Let $(E,D)$ be a rank-$\ell$ metric vector bundle with compatible connection and $W\ge0$ a smooth scalar potential. For the Dirichlet-exhaustion negative index of
\[
 \mathscr L=D^*D+q\,\mathrm{Id}_E-W\,\mathrm{Id}_E
\]
there is a dimensional constant $B_n$ such that
\begin{equation}\label{eq:CLR}
 N_-(\mathscr L)\le\ell\,a B_n\int_M W^{n/2}\,d\mu.
\end{equation}
The same estimate applies to $D^*D+q\,\mathrm{Id}_E+\mathscr V$, with $W$ replaced by $(-\lambda_{\min}(\mathscr V))_+$, for any smooth symmetric endomorphism $\mathscr V$.
\end{proposition}
\begin{proof}
First fix a relatively compact smooth domain $\Omega$ and Dirichlet boundary conditions. Kato's inequality and scalar semigroup domination show that
\[
 |e^{-t(D^*D+q)}s|\le e^{-t(\Delta_\Omega+q)}|s|,
\]
pointwise. The latter scalar heat kernel is bounded by the global minimal kernel of $\Delta+q$ by Dirichlet domain monotonicity. For completeness, we recall the quantitative part of the Birman--Schwinger/Lieb trace estimate underlying \cite[Theorem 37 and formula (26)]{BerardBesson}. Choose an integer $m>n/2-1$. Its trace regularization uses a nonnegative function $\Phi_m$ which is $O(s^{m+1})$ as $s\downarrow0$ and $O(s)$ as $s\to\infty$ (one may take the convex majorant of $s(1-e^{-s})^m$ constructed in the cited argument). Heat-semigroup domination gives, for a dimensional regularization factor $D_m$,
\[
 N_-\bigl(\mathscr L^D_\Omega\bigr)
 \le D_m\ell\int_\Omega\int_0^\infty
   p_{q,\Omega}(t,x,x)\,\Phi_m(tW(x))\,\frac{dt}{t}\,d\mu(x).
\]
The point of retaining the nonnegative scalar potential $q$ is that the domination and this same trace bound hold for the reference operator $D^*D+q$. Since $p_{q,\Omega}(t,x,x)\le at^{-n/2}$, a change of variables $s=tW(x)$ yields
\[
 N_-\bigl(\mathscr L^D_\Omega\bigr)
 \le D_m\ell a
 \left(\int_0^\infty s^{-1-n/2}\Phi_m(s)\,ds\right)
 \int_\Omega W^{n/2}\,d\mu.
\]
The integral is finite: at zero $m+1>n/2$, and at infinity $n>2$. Absorbing it into $B_n$ proves the claimed inequality on $\Omega$.
The kernel domination is unaffected by adding a nonnegative scalar $q$ on both sides, and the exponents in the trace argument depend only on $n$. Taking the supremum over Dirichlet exhaustions yields \eqref{eq:CLR}. For a general $\mathscr V$, the quadratic-form inequality $\mathscr V\ge -W\,\mathrm{Id}_E$ and the min--max principle give the final assertion. This is the reference-potential version of the classical CLR argument; see also \cite{MolchanovVainberg} for the role of diagonal heat-kernel estimates in general CLR theory.
\end{proof}

\begin{remark}[Spectral precision]\label{rem:clr}
The quantity $N_-$ counts \emph{strictly negative} Dirichlet eigenvalues. It does not automatically count zero eigenvalues. A separate compactly supported perturbation is required before applying \cref{prop:CLR} to $L^2$ harmonic forms; see \cref{lem:zero}.
\end{remark}

\section{\texorpdfstring{$L^2$-harmonic forms}{L2-harmonic forms} and the Gauss equation}
For $1\le p\le n-1$, the Hodge Laplacian on $p$-forms satisfies the Weitzenb\"ock formula
\begin{equation}\label{eq:weit}
 \Delta_p=\nabla^*\nabla+\mathscr R_p.
\end{equation}
On an immersed Euclidean submanifold the Gauss equation gives $|\mathrm{Rm}^M|\le C_n|A|^2$, so the algebraic Weitzenb\"ock operator satisfies
\begin{equation}\label{eq:Rbound}
 \|\mathscr R_p\|_{\mathrm{op}}\le d_{n,p}|A|^2.
\end{equation}
Here $d_{n,p}$ depends only on $n,p$; no constant depends on the ambient codimension.

\begin{lemma}[Counting zero modes by negative perturbation]\label{lem:zero}
If $\phi\in C_c^\infty(M)$ is nonnegative and nonzero, then for every $\varepsilon>0$,
\[
 \dim\cH_{(2)}^p(M)\le N_-(\Delta_p-\varepsilon\phi\,\mathrm{Id}),
\]
where the right-hand side is the Dirichlet-exhaustion negative index.
\end{lemma}
\begin{proof}
Let $V\subset\cH_{(2)}^p(M)$ be finite dimensional. By elliptic regularity, its elements are smooth. Unique continuation for the Hodge Laplacian (see \cite{AKS}) implies
\[
 c_V:=\inf_{\omega\in V,\ \|\omega\|_2=1}\int_M\phi|\omega|^2>0.
\]
Take compactly supported Lipschitz cutoffs $\chi_R(x)=\eta(d(o,x)/R)$, with $\chi_R=1$ near $\supp\phi$ and $|d\chi_R|\le2/R$. Completeness makes their supports compact. For harmonic $\omega$, the product rules give $d(\chi_R\omega)=d\chi_R\wedge\omega$ and $d^*(\chi_R\omega)=-\iota_{\nabla\chi_R}\omega$; the Clifford identity yields
\[
 |d\chi_R\wedge\omega|^2+|\iota_{\nabla\chi_R}\omega|^2
 =|d\chi_R|^2|\omega|^2.
\]
Consequently, for the quadratic form of $\Delta_p-\varepsilon\phi$,
\[
 Q_\varepsilon[\chi_R\omega]
 \le \frac4{R^2}\|\omega\|_2^2-\varepsilon c_V\|\omega\|_2^2.
\]
For large $R$ this is negative definite on $\chi_RV$, and $\omega\mapsto\chi_R\omega$ is injective by unique continuation. Smooth approximation permits the same conclusion in the Dirichlet form domain. Since $V$ was arbitrary, the result follows.
\end{proof}

\begin{proposition}[Quantitative $L^2$ cohomology]\label{prop:L2}
For $1\le p\le n-1$,
\[
 \dim\cH_{(2)}^p(M)\le C_{n,p}\int_M|A|^n\,d\mu.
\]
\end{proposition}
\begin{proof}
Use $q=|H|^2$ and rewrite
\[
 \Delta_p-\varepsilon\phi
 =\nabla^*\nabla+q\,\mathrm{Id}+(\mathscr R_p-q\,\mathrm{Id}-\varepsilon\phi\,\mathrm{Id}).
\]
By \eqref{eq:Rbound} and $|H|^2\le n|A|^2$, the negative part of the least eigenvalue of the bracketed endomorphism is bounded by
\[
 W_\varepsilon=(d_{n,p}+n)|A|^2+\varepsilon\phi.
\]
Apply \cref{prop:heat,prop:CLR,lem:zero} to obtain
\[
 \dim\cH_{(2)}^p(M)
 \le \binom{n}{p}\,a_nB_n\int_M W_\varepsilon^{n/2}\,d\mu.
\]
For $0<\varepsilon\le1$, the integrand is dominated by $C_{n,p}(|A|^n+\phi^{n/2})\in L^1(M)$. Letting $\varepsilon\downarrow0$ proves the assertion, with for instance
\[
 C_{n,p}=\binom{n}{p}a_nB_n(d_{n,p}+n)^{n/2}.
\]
\end{proof}

\section{Ends: the potential-theoretic input}
A nonparabolic end means an end whose equilibrium potential is nonconstant, equivalently an end of positive capacity at infinity. This is not the same as positivity of the capacity of a compact set with Dirichlet data on the \emph{interior} boundary of an arbitrary end domain: that weaker assertion does not characterize nonparabolicity.

\begin{lemma}[Nonparabolic ends under finite total mean curvature]\label{lem:np}
If $n\ge3$ and $F:M^n\to\R^N$ is complete with $\int_M|H|^n<\infty$, then every end of $M$ is nonparabolic.
\end{lemma}
\begin{proof}
For orientable immersions this is the nonparabolicity theorem of Fu--Li \cite{FuLi}. The proof there uses Michael--Simon inequalities on the tails and the appropriate end-capacity criterion; in particular it does not identify relative Dirichlet capacity at the inner boundary with capacity at infinity. If $M$ is nonorientable, pull the immersion back to its orientable double cover. The cover is complete and has twice the finite total mean-curvature integral. Each end of the double cover is therefore nonparabolic. To see that this descends, note that if an end downstairs were parabolic, its equilibrium-potential cutoff sequence could be lifted to any component of its finite-sheeted preimage, with Dirichlet energies multiplied by at most two. The lifted end would then be parabolic as well, a contradiction.
\end{proof}

\begin{proof}[Proof of \cref{thm:general}]
By $|H|\le\sqrt n|A|$, finite $\TT(F)$ implies $\int|H|^n<\infty$; thus \cref{lem:np} makes all ends nonparabolic. Li--Tam theory \cite{LiTam} gives
\[
 e(M)=e_\np(M)\le\dim\cH_D^0(M),
\]
where $\cH_D^0(M)$ is the space of bounded harmonic functions of finite Dirichlet energy. Differentiation maps this space into $\cH_{(2)}^1(M)$ with kernel precisely the constants. Hence
\[
 e(M)\le1+\dim\cH_{(2)}^1(M)\le1+C_{n,1}\TT(F),
\]
using \cref{prop:L2}. The remaining $L^2$ bounds follow from \cref{prop:L2}. Since $e(M)\ge1$ for a connected noncompact manifold, small total curvature forces one end; integrality of the dimensions gives the stated $L^2$ gap as well.
\end{proof}

\section{Geometric preliminaries for the Morse method}
\begin{definition}[Asymptotically planar graphical ends]\label{def:ends}
An immersion $F:M^n\to\R^N$ has $e$ \emph{asymptotically planar graphical ends} if there is a compact connected smooth codimension-zero submanifold $K\Subset M$ with boundary, such that
\[
M\setminus\operatorname{int}K\cong \bigsqcup_{a=1}^{e}\bigl(\SSS^{n-1}\times[R_a,\infty)\bigr),
\]
and, after a rigid ambient motion depending on $a$, the immersion of the $a$th end is given by
\begin{equation}\label{eq:graph}
F_a(x)=(x,Y_a(x))\in\Pi_a\oplus\Pi_a^\perp,
\qquad |x|\ge R_a,
\end{equation}
where $\Pi_a$ is an $n$-plane, $Y_a\in C^2(\{|x|\ge R_a\};\Pi_a^\perp)$, and
\begin{equation}\label{eq:graphdecay}
\lim_{r\to\infty}\sup_{|x|\ge r}
\left(\frac{|Y_a(x)|}{|x|}+|DY_a(x)|+|x|\,|D^2Y_a(x)|\right)=0.
\end{equation}
We allow different tangent planes $\Pi_a$ and intersections of the \emph{images} of distinct ends. The end charts are embeddings individually.
\end{definition}

This hypothesis implies that the immersion is proper: on each end, $|F_a(x)|\ge|x|\to\infty$, and the remaining part of $M$ is compact. It also implies that the inclusion $K_R\hookrightarrow M$ is a homotopy equivalence, where $K_R$ is obtained by truncating each graph at $|x|=R$ for any sufficiently large common $R$.

Put $k=N-n\ge1$. Let $U_{\nu} M$ be the unit normal sphere bundle and define the unit normal Gauss map
\[
\Gamma:U_{\nu} M\longrightarrow\SSS^{N-1},\qquad
\Gamma(p,\xi)=\xi.
\]
For $\xi\in\nu_pM$ of unit norm let $S_\xi:T_pM\to T_pM$ be the shape operator, defined by
$\langle S_\xi X,Y\rangle=\langle A(X,Y),\xi\rangle$.
\begin{definition}
Following Chern--Lashof \cite{ChernLashof}, the normalized total absolute curvature is
\begin{equation}\label{eq:TA}
\TA(F):=\frac1{|\SSS^{N-1}|}
\int_{U\nu M}|\det S_\xi|\,d\xi\,d\mu(p).
\end{equation}
\end{definition}

\begin{lemma}[Total-curvature domination]\label{lem:TAcurvature}
For any smooth immersed $n$-manifold in $\R^N$ with $N>n$,
\begin{equation}\label{eq:TAdom}
\TA(F)\le c_{n,N}\int_M|A|^n\,d\mu,
\end{equation}
with $c_{n,N}$ as in \eqref{eq:cconst}.
\end{lemma}
\begin{proof}
If $\lambda_1,\dots,\lambda_n$ are the eigenvalues of the symmetric operator $S_\xi$, the arithmetic--geometric mean inequality gives
\[
|\det S_\xi|=\prod_i|\lambda_i|
\le\left(\frac{\sum_i\lambda_i^2}{n}\right)^{n/2}
=n^{-n/2}|S_\xi|^n\le n^{-n/2}|A|^n.
\]
Integrating over each normal fiber $\SSS^{k-1}$ and dividing by $|\SSS^{N-1}|$ proves the assertion.
\end{proof}

\section{Generic height functions and compact truncations}
For $v\in\SSS^{N-1}$ consider the height function
\[
h_v(p):=\langle F(p),v\rangle.
\]
Write $v_a=P_{\Pi_a}v$. Let $\mathscr G\subset\SSS^{N-1}$ be the set of $v$ such that (i) $v_a\ne0$ for every end, and (ii) $v$ is a regular value of $\Gamma$. Since there are finitely many ends, condition (i) excludes finitely many lower-dimensional great subspheres; condition (ii) fails on a measure-zero set by Sard's theorem. Hence $\mathscr G$ has full measure.

\begin{lemma}[Interior critical points]\label{lem:crit}
For every $v\in\mathscr G$, the height function $h_v$ has finitely many critical points, all nondegenerate, and all contained in a compact subset of $M$.
\end{lemma}
\begin{proof}
On the $a$th end, let $x\in\Pi_a$ and write $v=v_a+v_a^\perp$. Then
\[
d(h_v\circ F_a)_x=v_a+(DY_a(x))^*v_a^\perp.
\]
By \eqref{eq:graphdecay}, for all sufficiently large $|x|$ (depending on $v$),
$|d(h_v\circ F_a)_x|\ge |v_a|/2>0$. Thus all critical points lie in a compact subset. At a critical point $p$, the vector $v$ is normal and
\[
\operatorname{Hess}_M h_v(X,Y)=\langle A(X,Y),v\rangle=\langle S_vX,Y\rangle.
\]
Since $v$ is a regular value of the unit normal Gauss map, $\det S_v\ne0$ at every preimage of $v$. Consequently all critical points are nondegenerate. A compact discrete set of critical points is finite.
\end{proof}

For large $R$, write $K_R$ for the compact core truncated along the $e$ sections
\[
\Sigma_{a,R}=\{(R\theta,Y_a(R\theta)):\theta\in\SSS^{n-1}\}.
\]
The restriction of $h_v/R$ to the boundary component $\Sigma_{a,R}$, in the $\theta$ coordinates, is
\begin{equation}\label{eq:boundaryg}
g_{a,R}(\theta)
=\langle\theta,v_a\rangle+R^{-1}\langle Y_a(R\theta),v_a^\perp\rangle.
\end{equation}

\begin{lemma}[Two boundary critical points]\label{lem:boundary}
Fix $v\in\mathscr G$. For sufficiently large $R$, the restriction $h_v|_{\Sigma_{a,R}}$ has exactly two nondegenerate critical points $p_{a,-}$ and $p_{a,+}$, of boundary Morse indices $0$ and $n-1$, respectively. If $\nu_R$ is the outward unit conormal of $K_R$, then
\[
\partial_{\nu_R}h_v(p_{a,-})<0,
\qquad
\partial_{\nu_R}h_v(p_{a,+})>0.
\]
In particular, each end contributes exactly one \emph{type-N} boundary critical point, and its boundary Morse index is zero. The same sufficiently large truncation radius can be chosen simultaneously for all ends (for the fixed $v$).
\end{lemma}
\begin{proof}
By the chain rule and \eqref{eq:graphdecay}, the error term in \eqref{eq:boundaryg} tends to zero in $C^2(\SSS^{n-1})$. Thus $g_{a,R}$ converges in $C^2$ to the nonzero linear height function $\theta\mapsto\langle\theta,v_a\rangle$. The latter has precisely two nondegenerate critical points, at $-v_a/|v_a|$ (minimum, index $0$) and $v_a/|v_a|$ (maximum, index $n-1$). Stability of nondegenerate critical points and a lower bound for the gradient away from their neighborhoods give the same critical-point pattern for $g_{a,R}$ when $R$ is large.

The outward conormal on the graph is $\nu_R=(\theta,0)+o(1)$, uniformly as $R\to\infty$ in the ambient Euclidean norm. This follows by solving the orthogonality equations for the tangent vectors $(X,DY_a(X))$ and normalizing the radial tangent direction. The error is controlled by $\sup_{|x|=R}|DY_a(x)|$. Therefore
$\partial_{\nu_R}h_v=\langle v_a,\theta\rangle+o(1)$.
Its value tends to $-|v_a|$ at the minimum and $+|v_a|$ at the maximum, giving the signs claimed. Because the number of ends is finite, taking the maximum of the finitely many thresholds yields a common admissible $R(v)$.
\end{proof}

\section{Absolute boundary Morse homology and a rank estimate}
We recall the relevant standard consequence of Morse theory for compact manifolds with boundary, in the form proved by Laudenbach \cite{Laudenbach}; for the classical closed case see \cite{Milnor}.

\begin{lemma}[Absolute Morse counting]\label{lem:Morse}
Let $K$ be a compact connected smooth $n$-manifold with boundary and let $f\in C^\infty(K)$ have nondegenerate interior critical points. Suppose $f|_{\partial K}$ is Morse and $\partial_\nu f\ne0$ at all its critical points, with $\nu$ the outward conormal. For every field $\F$ there is a finite chain complex computing $H_*(K;\F)$ whose generators in degree $j$ consist of
\begin{itemize}[leftmargin=2em]
\item interior critical points of $f$ of Morse index $j$;
\item boundary critical points of $f|_{\partial K}$ of index $j$ with $\partial_\nu f<0$.
\end{itemize}
\end{lemma}
\begin{proof}[Justification of the coefficient-field statement]
Laudenbach's Theorem A constructs a complex $C_*^{N}(f;\mathbb Z)$ of free abelian groups with the indicated generators and homology $H_*(K;\mathbb Z)$ \cite{Laudenbach}. Tensoring with an arbitrary field $\F$ produces the complex $C_*^{N}(f;\mathbb Z)\otimes_{\mathbb Z}\F$. The universal coefficient theorem for chain complexes of free abelian groups identifies its homology with $H_*(K;\F)$ (both sides have the same homology and torsion contributions). Thus neither orientability nor a characteristic-zero restriction is needed. The auxiliary Morse--Smale pseudo-gradient of Laudenbach may be chosen after fixing $f$ and does not change the set or indices of its critical points.
\end{proof}
\begin{remark}
The sign in this statement is essential. For example, for a closed ball and a nonconstant linear function there are no interior critical points; its boundary minimum is the single generator for $H_0$. Its outward normal derivative is negative. The boundary maximum has positive outward normal derivative and is not included in the absolute complex.
\end{remark}

\begin{proposition}[A pointwise Morse inequality with an end correction]\label{prop:pointwise}
For any $v\in\mathscr G$, let $\mu(v)$ be the number of critical points of $h_v$ on $M$ and let $m_0(v)$ be the number of those with Morse index $0$. Then, for every field $\F$,
\begin{equation}\label{eq:pointwise}
\bb(M;\F)+e-2+2m_0(v)\le\mu(v).
\end{equation}
In particular, $\bb(M;\F)+e-2\le\mu(v)$.
\end{proposition}
\begin{proof}
Choose $R$ large enough for \cref{lem:crit,lem:boundary}. Coincident critical values can be separated by an arbitrarily small perturbation supported near the finite set of interior critical points, without changing their Morse indices or the boundary critical-point pattern. The Morse--Smale pseudo-gradient is then chosen for this fixed function; the same generator counts result. This also avoids imposing a distinct-critical-value condition on the original height function. The inclusion $K_R\hookrightarrow M$ is a homotopy equivalence, so $\bb(K_R;\F)=\bb(M;\F)$. By \cref{lem:Morse}, the absolute Morse chain complex $(C_*,\partial_*)$ of $h_v|_{K_R}$ has exactly $\mu(v)$ interior generators and $e$ boundary generators. By \cref{lem:boundary}, every boundary generator has degree $0$, and there is exactly one per end. In particular,
\begin{equation}\label{eq:chaincounts}
\sum_j\dim C_j=\mu(v)+e,
\qquad\dim C_0=m_0(v)+e.
\end{equation}
Because $K_R$ is connected, $\dim H_0(C_*)=1$. Thus
\begin{equation}\label{eq:rankdiff}
\rank\partial_1=\dim C_0-1=m_0(v)+e-1.
\end{equation}
For a finite complex of vector spaces over $\F$, the rank--nullity identity gives
\[
\sum_j\dim H_j(C_*)
=\sum_j\dim C_j-2\sum_{j\ge1}\rank\partial_j.
\]
Using \eqref{eq:chaincounts} and \eqref{eq:rankdiff},
\begin{align*}
\bb(M;\F)
&=\mu(v)+e-2\sum_{j\ge1}\rank\partial_j\\
&\le\mu(v)+e-2\bigl(m_0(v)+e-1\bigr)\\
&=\mu(v)-e+2-2m_0(v),
\end{align*}
which is \eqref{eq:pointwise}.
\end{proof}

The use of the \emph{strong} chain-complex identity rather than just the weak Morse inequalities accounts for the improvement from $\bb(M)\le\TA(F)+e$ to $\bb(M)+e-2\le\TA(F)$.

\section{The Gauss-map area formula and the main theorem}
\begin{lemma}[Counting by total absolute curvature]\label{lem:average}
For the immersion above,
\begin{equation}\label{eq:average}
\TA(F)=\frac1{|\SSS^{N-1}|}
\int_{\SSS^{N-1}}\mu(v)\,dv,
\end{equation}
where the counting function is defined almost everywhere by the number of points in the inverse image of the unit normal Gauss map (equivalently, the number of critical points of the ambient height function).
\end{lemma}
\begin{proof}
The domain and codomain of $\Gamma:U\nu M\to\SSS^{N-1}$ both have dimension $N-1$. Its absolute Jacobian is $|\det S_\xi|$, as follows from the Weingarten formula in an orthonormal splitting of horizontal and vertical spaces. The area formula \cite{Federer} applies to the smooth map $\Gamma$ even when its domain is noncompact: apply the formula to a countable exhaustion and then use monotone convergence, since the integrands are nonnegative. It yields
\[
\int_{U\nu M}|\det S_\xi|\,d\xi\,d\mu
=\int_{\SSS^{N-1}}\#\Gamma^{-1}(v)\,dv.
\]
For almost every $v$, the preimages of $v$ are in bijection with critical points of $h_v$, since $dh_v(p)=0$ precisely when $v$ is normal to $T_pM$. This proves the formula.
\end{proof}

\begin{proof}[Proof of \cref{thm:ordinary}]
By \cref{prop:pointwise},
$\bb(M;\F)+e-2\le\mu(v)$ for almost every $v\in\SSS^{N-1}$. Averaging and applying \cref{lem:average} gives
\[
\bb(M;\F)+e-2\le\TA(F).
\]
To deduce the individual bounds, note that for any connected compact $n$-manifold $K_R$ with $e$ boundary components diffeomorphic to $\SSS^{n-1}$ and $n\ge2$,
\begin{equation}\label{eq:topineq}
b_{n-1}(K_R;\F)\ge e-1.
\end{equation}
Indeed the homology exact sequence of $(K_R,\partial K_R)$ contains
\[
H_n(K_R,\partial K_R;\F)\longrightarrow
H_{n-1}(\partial K_R;\F)
\longrightarrow H_{n-1}(K_R;\F),
\]
where the middle space has dimension $e$ and the left space has dimension at most $1$: for connected $K_R$, Poincar\'e--Lefschetz duality with the orientation local system identifies it with $H^0(K_R;\mathcal O\otimes\F)$, which has dimension either $0$ or $1$. This argument works over any field and for nonorientable $K_R$. Thus the image in $H_{n-1}(K_R;\F)$ has dimension at least $e-1$. Since $K_R$ is connected and $n\ge2$, we have $b_0(K_R)=1$ in a degree different from $n-1$. Consequently,
\begin{equation}\label{eq:be}
\bb(M;\F)\ge e.
\end{equation}
Combining \eqref{eq:be} with \eqref{eq:morsecl}, we obtain
$2e-2\le\TA(F)$, hence
\[
e\le1+\frac12\TA(F).
\]
As $e\ge1$, \eqref{eq:morsecl} also implies
$\bb(M;\F)\le1+\TA(F)$. Finally apply \cref{lem:TAcurvature} to obtain \eqref{eq:ordinarybound} and \eqref{eq:planarend}.
\end{proof}

\begin{remark}[Normalization in dimension two]
For the classical catenoid in $\R^3$, viewed as an immersed cylinder, $b_0=b_1=1$, $e=2$ and $\int_M|A|^2=8\pi$. Since its total absolute Gaussian curvature is $4\pi$ and the unit normal fiber consists of two points, $\TA(F)=2$. Thus
\[
 \bb(M;\F)+e-2=2=\TA(F)=c_{2,3}\int_M|A|^2,
 \qquad c_{2,3}=\frac1{4\pi}.
\]
This illustrates the normalization of total absolute curvature and equality in the refined Morse bound. It does not invoke the higher-dimensional minimal-end regularity theorem.
\end{remark}

\section{Minimal immersions and a weighted-end calculation}\label{sec:minimal}
We isolate the derivative step in the application of \cite[Theorem 3.2]{DingZhang}. It is important that a pointwise $O(r^{-n})$ remainder cannot be differentiated without additional regularity: the derivative estimates used below come instead from the weighted $C^2$ norm in the proof of their theorem.

\begin{lemma}[From cylindrical $C^2$ decay to graphical decay]\label{lem:weights}
Let $Y:\{x\in\R^n:|x|>R_0\}\to\R^k$ be $C^2$ and set
\[
 t=\log r,\quad \theta=x/|x|,\quad
 U(t,\theta)=e^{-t}Y(e^t\theta).
\]
If $\|U\|_{C^2(\SSS^{n-1}\times[t,t+1])}\to0$ as $t\to\infty$, then
\[
 \lim_{r\to\infty}\sup_{|x|\ge r}
 \left(\frac{|Y(x)|}{|x|}+|DY(x)|+|x|\,|D^2Y(x)|\right)=0.
\]
More quantitatively, if the cylindrical $C^2$ norm is $O(e^{-\gamma t})$ for some $\gamma>0$, the displayed expression is $O(r^{-\gamma})$.
\end{lemma}
\begin{proof}
For a scalar or vector-valued function $V(t,\theta)$, differentiation of $V(\log r,x/r)$ in an Euclidean coordinate direction is a linear combination of $r^{-1}\partial_tV$ and $r^{-1}\nabla_{\SSS^{n-1}}V$, with coefficients bounded uniformly on the unit sphere. Apply this to $Y(r\theta)=rU(\log r,\theta)$. Differentiating once gives, in Euclidean tensor norms,
\[
 |DY(r\theta)|\le C_n\bigl(|U|+|\partial_tU|+|\nabla_{\SSS^{n-1}}U|\bigr).
\]
A second Euclidean derivative produces a factor $r^{-1}$ and involves at most two cylindrical derivatives of $U$, so
\[
 r|D^2Y(r\theta)|\le C_n\sum_{j=0}^2|\nabla_{\rm cyl}^jU(t,\theta)|.
\]
Since $|Y(r\theta)|/r=|U(t,\theta)|$, all three quantities are bounded by $C_n\|U\|_{C^2(\SSS^{n-1}\times[t,t+1])}$. Uniform decay on successive unit cylinders implies uniform decay on every tail. The exponential rate yields the stated polynomial rate because $e^{-\gamma t}=r^{-\gamma}$.
\end{proof}

\begin{proof}[Proof of Corollary~\ref{cor:minimal}]
For $n\ge3$, Ding--Zhang's regular-at-infinity theorem \cite[Theorem 3.2]{DingZhang}, building on Anderson \cite{Anderson}, asserts that a complete connected minimal immersion with $\int|A|^n<\infty$ has a finite collection of single-sheeted graphical ends over exterior domains in affine $n$-planes and is properly immersed. More precisely, the proof of the implication from planar tangent cones to regularity in \cite[Section 3]{DingZhang} first obtains the cylindrical end parametrization
\[
 (t,\theta)\longmapsto(e^t\theta,e^tU(t,\theta))
\]
with $\|U\|_{2,\gamma}<\infty$ for some $\gamma>0$; in their notation this is the weighted $C^2$ norm. Hence
\[
 \|U\|_{C^2(\SSS^{n-1}\times[t,t+1])}=O(e^{-\gamma t}).
\]
Lemma~\ref{lem:weights} now verifies condition \eqref{eq:graphdecay} for each end. The affine displacement of its asymptotic plane adds only a term $O(r^{-1})$ to $U$ and does not affect this conclusion. Taking a common large truncation gives the compact core of Definition~\ref{def:ends}. \cref{thm:ordinary} gives the two estimates stated in \cref{cor:minimal}; properness and finite topological type follow from the same graphical-end structure.
\end{proof}

\begin{proof}[Proof of Corollary~\ref{cor:hsl}]
Identify $\C^n$ with $\R^{2n}$ and apply \cref{thm:ordinary}. The Lagrangian and Hamiltonian stationary equations are not used in the Morse argument.
\end{proof}

\begin{corollary}[A curvature threshold for multiple ends]\label{cor:threshold}
Under the hypotheses of \cref{thm:ordinary}, if $e(M)\ge r\ge2$ and $\TT(F)<\infty$, then
\[
 \TT(F)\ge\frac{2(r-1)}{c_{n,N}}.
\]
In particular, $\TT(F)<2/c_{n,N}$ forces exactly one end.
\end{corollary}
\begin{proof}
Rearrange \eqref{eq:planarend}.
\end{proof}

\section{Comparison of the two methods}\label{sec:comparison}
\subsection{Distinct topological information}
The spectral route estimates the dimension of the kernel of a Hodge Laplacian in $L^2$; it gives bounds for the \emph{reduced $L^2$-Betti numbers}. Without additional information about the ends, those quantities cannot be substituted for ordinary Betti numbers. Its end count is mediated by the Li--Tam inequality
\[
 e(M)=e_{\mathrm{np}}(M)\le1+b^1_{(2)}(M),
\]
and it remains applicable even for a nonproper immersion. The argument controls the negative part of the Weitzenb\"ock curvature term through the Gauss equation and is insensitive to the ambient codimension except through the geometric immersion itself.

The Morse route, by contrast, counts the critical points of \emph{ambient height functions}. A compact truncation is available only after the ends are identified geometrically. Each graph end contributes one degree-zero boundary generator to the absolute Morse complex; connectedness forces a rank $e-1$ condition on the first differential. The resulting pointwise inequality
\[
 \bb(M;\F)+e(M)-2+2m_0(v)\le\mu(v)
\]
contains strictly more ordinary-topological information than a $\dim\cH^1_{(2)}$ estimate. Averaging $\mu(v)$ over directions is a geometric coarea/area-formula step, rather than an elliptic spectral estimate.

\subsection{Assumptions, constants, and scope}
The heat-kernel argument requires $n\ge3$ and $\TT(F)<\infty$ but no prescribed geometry at infinity. Its dimensional coefficients arise from Michael--Simon and the abstract CLR estimate; they are not explicitly optimized here. The Morse method needs a finite family of single-sheeted graph ends satisfying the scale-invariant $C^2$ decay in \cref{def:ends}. Within that class its constants are explicit and the additional conclusions concern ordinary homology. No numerical ordering between the two end-count constants is asserted, because the spectral constant is not computed explicitly.

The distinction is particularly relevant for nonminimal immersions. Finite $\TT(F)$ gives the end and $L^2$ bounds of \cref{thm:general} but does not, from the results used here alone, imply the graphical-end hypothesis of \cref{thm:ordinary}. The curvature-tameness condition of Bessa--Jorge--Montenegro \cite{BJM} is one route to properness and finite topological type, but is not equivalent to the spherical single-sheeted condition. For minimal immersions, regularity \cite{Anderson,DingZhang} makes both routes available and the Morse route provides the stronger explicit conclusion.

\subsection{The critical dimension and possible extensions}
The spectral bound relies on $n/2>1$, the same exponent restriction as the classical CLR inequality. The present article accordingly assumes $n\ge3$. The geometric Morse argument under \emph{given} graphical-end hypotheses can also be formulated in dimension two, but finite total curvature of a minimal surface does not automatically give a single-valued graph end in that dimension. A higher-dimensional finite-curvature assumption alone has not been shown here to yield uniform ordinary Betti-number estimates for arbitrary nonminimal immersions.
\subsection{A scale-invariant gap observation}
The Morse bound implies that a planar-ended immersion with $e(M)\ge r\ge2$ satisfies
\[
 \TT(F)\ge\frac{2(r-1)}{c_{n,N}},
\]
provided $\TT(F)<\infty$. In the general setting, the spectral inequality only gives the dimensional (nonexplicit) threshold $\TT(F)\ge(r-1)/C_n$. Both inequalities are invariant under $F\mapsto\lambda F$. These are statements about necessary curvature cost, not optimality claims.

\end{document}